# Constructing Piecewise Flat Pseudo-Manifolds with Minimal Pseudo-Foliations

Simon P Morgan

University of Minnesota

**Abstract**

Piecewise Euclidean structures (identified solid Euclidean polyhedra) on topological 3-dimensional manifolds and pseudo-manifolds are constructed so that they admit pseudo-foliations, a generalized type of foliation. The construction of non-manifold point neighborhoods is done to preserve as much of the geometric, and vector space, structure of Euclidean space as possible. This enables structures, such as foliations and calibrations, to generalize over to pseudo-manifolds.

The main result is that the cone of any compact topological surfaces can be given a piecewise flat metric structure that admits a pseudo-foliation by minimal surfaces. In some situations orientation reversing holonomy is an obstruction and in others it is just the opposite, a sufficient condition. Also foliations by minimal surfaces can extend across connect sum operations.

## 1  Introduction

The main examples of this paper are cones of surfaces. These are generally pseudo-manifolds, as the link of the cone point need not be a sphere. Such spaces can occur naturally as orbifolds and as one point compactifications of knot compliments. For example, the cone of $\mathbb{RP}^2$ is the quotient of $\mathbb{R}^3$ by the involution $x \sim -x$.

The approach of this paper comes from the following observation relating different approaches. The technique of calibration [Harvey and Lawson] can be used to determine a foliation by minimal surfaces. [Morgan] gives a brief review of novel uses of calibrations. Foliations have also been studied with a more topological approach. Early examples include [Thurston], [Sullivan], [Banchoff] and more recently the work of David Gabai [Gabai],[Gabai and Kazez]. There is a standard condition that the foliation must be orthogonal to any boundary. However once we have this way of dealing with foliations, we can go one step further and identify boundaries to construct new foliated spaces from old, subject to holonomy conditions. This type of construction of piecewise flat structures, not particularly associated with foliations was carried out by [Jones].

This work has similarities to work with surfaces in 3 manifold theory. Embedded surfaces play an important role in 3-manifold theory as they can decompose a complicated 3-manifolds into simpler components. For this to work the surfaces need to be as simple as possible, to give decompositions that are as simple as possible. [Meeks, Simon and Yau] use minimal surface theory to find suitable embedded surfaces. A combinatorial approach uses normal surface theory. These are surfaces that intersect triangulating tetrahedra in canonical ways, e.g.: [Jaco and Rubinstein] and [Stocking]. The work in this paper has things in common with both approaches, and has the potential to allow isotopy of minimal surfaces in pseudo-manifolds to be studied. This isotopy has the unusual property of allowing the topology of a surface to change as it passes over a non-manifold point. See the remark in section 5.

### 1.1 Preliminaries

It is worth mentioning two preliminary issues. The first is topological obstruction, such as the Poincaré-Hopf index theorem. That is a boundary-less surface cannot admit a nowhere vanishing continuous vector field unless its Euler characteristic is zero. However if we add boundary to the surface with no condition on how the vector field meets the boundary, the obstruction vanishes. The second issue is singular local curvature. Consider a cone point in an otherwise flat surface. For what values of cone angle is there a foliation by geodesics? Only when the cone angle is $2\pi$ can a foliation by geodesics exist on the cone. As discuss below for values of $n\pi$ a singular foliation can be admitted with a branch point of order n for the geodesics going through the apex. We will say that a geodesic terminating at an apex is a branch point of order 1. A branch point of order 2 is a regular point. The order of the branch point is analogous to the index of a vector field in the Poincaré-Hopf index theorem, and by the Gauss-Bonnet theorem we know that curvature is connected to Euler Characteristic. So for a piecewise flat surface with cone points of angle $n\pi$ only we have.

$$\text{Euler characteristic} = \text{Total angle deficit at vertices}$$

$$= \Sigma\ 2\text{-}n \text{ over branch points of foliation.} \qquad (1)$$

Hence global topology and local geometry are related and can act as obstructions to minimal foliations. [Thurston] mentions this in purely topological terms.

### 1.2 Holonomy for minimal surfaces and calibration for area minimization.

The version of holonomy we shall use here is to take a basis of vectors at a point in the Euclidean space and parallel transport them through the space around loops. One cycle around a loop will have put the basis through a rigid motion, an element of $O(n)$. The existence of a fixed point in the action of the element of $O(n)$ on $\mathbb{S}^{n-1}$ indicates the possibility of a foliation by minimal surfaces in a neighborhood of that loop. The fixed point corresponds to a normal direction of the foliation.

Calibration is a global structure such that the existence of a calibration form indicates that a leaf of the foliation is not only minimal, but globally area minimizing. Globally area minimizing means that for any simple closed curve in the leaf, among the surfaces that have the closed curve as a boundary and have minimal area, one is contained in the leaf.

**1.3 The construction by surgery method**

Many of the proofs rely on a method of construction by surgery whereby a set S, of codimension 1 with boundary is removed. See figures 1, 8, 9 and 10. The neighborhood of the complement of this set is then compactified with two disjoint copies $S_1$ and $S_2$ of the original removed set. An identification pattern is then placed on $S_1$ and $S_2$ to complete the construction. We give an example in figure 1 and comment on the properties common with the constructions in figures 8, 9 and 10.

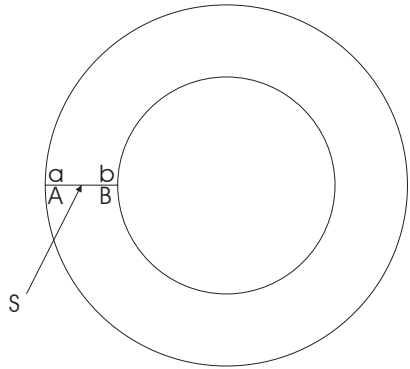

**Figure 1: Construction by surgery**

The annulus embedded in $\mathbb{R}^2$ in figure 1 contains the line segment S. If S is removed and its compliment compactified, $S_1$ compactifies the lower half with vertices A and B, while $S_2$ compactifies the top half with vertices a and b. $S_1$ and $S_2$ are disjoint, but if we tried to embed them with the same embedding as annulus, we could not, as they both fit into the same space where S was. So they will be isometrically immersed in $\mathbb{R}^2$ as a disjoint set. We can identify $S_1$ and $S_2$ to give us back the embedded annulus, a~A b~B, or we can obtain an embedded moebuis band, b~A and a~B.

**1.4 The paper outline**

Section 2 gives a range of elementary examples of constructed spaces and the spaces of foliations on them, including the implications of whether or not holonomy is orientation preserving or reversing. We give a construction to find we can foliate the hyperbolic compliment of the Borromean rings and other knots and links. Section 3 gives the formal framework for discussion of foliations on non-manifold points and even foliations which may have non-manifold points in the leaves. This generalizes the definition of foliation. Section 4 shows how to allocate metrics and foliate $\mathbb{S}^1$ bundles of surfaces and connect sums. In section 5 we give a generic construction for the cone of any topological surface so that it admits foliations by minimal surfaces. Finally we show that these surfaces are area minimizing using an adaptation of calibration proofs via the 'short cut lemma'.

# 2 Elementary examples of foliated surfaces and 3-manifolds.

**2.1 Orientation reversing holonomy: both obstruction and sufficient condition.**

Fixed point theorems for spheres ensure the existence of a fixed point in the orientation when the holonomy (as in 1.2) is orientation reversing. However in the orientation preserving case if the holonomy is ± the identity then all foliation directions are possible, otherwise there may be just one or no possible directions, depending on the dimension.

We can state the theorem:

**Theorem 1:** *A non-orientable n-disc bundle over $\mathbb{S}^1$ with a locally flat metric will admit at least one codimension one foliation by minimal submanifolds.*

**Proof:** The degree of the holonomy $h: \mathbb{S}^n \to \mathbb{S}^n$ is -1. and so it must carry one of its points, x, to its antipode; e.g. [Munkres] Thm 21.5. Let the normal vectors to the foliation point in directions x and –x. This set is preserved under the holonomy map. QED

Example, n=1. Any moebius band that is locally flat will admit a foliation by locally parallel lines. An easy counter example shows this to be false for the annulus. In each case we place identifications of the trapezoid shown in figure 2.

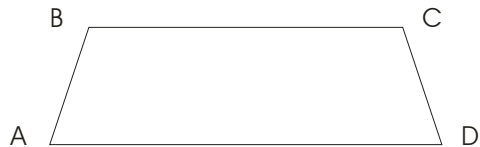

Figure 2

Consider first the non-orientable case. With the moebius band identifications A=C and B=D, both horizontal and vertical lines, but not other directions, will foliate the band. However with annular identifications B=C and A=D, the orientable case, no foliation will work. Although such an annulus is locally flat, it can be embedded in a cone which clearly cannot be foliated by geodesics that wind around the apex. The holonomy is a rotation of $180° - (\angle A + \angle D)$. In the special case where $\angle A + \angle D = 180°$, all directions are admissible for foliations.

For the orientable reversing holonomy, the theorem proves that whatever the dimensions of the flat trapezoid, a foliated moebius band can always be formed. Note we are not requiring any boundary conditions to be met by the foliation. Our example shows that orientation preserving holonomy can lead to no foliations. So in this context orientation reversing holonomy is a sufficient condition for one foliation to exist. However in special cases of orientation preserving holonomy, all directions can foliate. Orientation reversing holonomy acts as an obstruction to all directions foliating.

## 2.2 Foliating the hyperbolic Borromean rings compliment.

This construction, taken from the 'Not Knot' booklet [Epstein and Gunn], can be seen by examining the cube whose sides can be identified to give the Borromean rings complement. See figure 3, left.

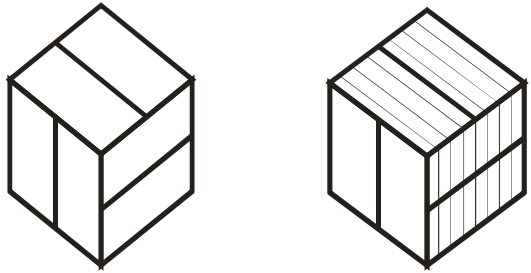

**Figure 3: Borromean rings compliment and its foliation**

To make the Borromean rings compliment, reflect each half face in that face's bisector. This identifies ends of bisectors thus forming the Borromean rings. Each ring is also an arc of singular positive curvature. To foliate this cube with identifications by minimal surfaces chose a set of surfaces parallel to one of the faces. See figure 3, right. After identification, each leaf is a minimal surface and there are two singular leaves which are discs whose boundaries are each one of the rings. These are the leaves which bisect the cube as shown in figure 3, right, and the leaf which forms 2 faces, front left, and right rear.

If, prior to identification the face bisectors are given a dihedral angle (see figure 4 left) and their lengths are reduced then eventually we achieve a rhombic dodecahedron (see figure 4 right). The line AA on the left in figure 4 becomes a point A on the sphere at infinity on the ideal rhombic dodecahedron on the right. These are explained more fully in the 'Not Knot' booklet and video [Epstein and Gunn].

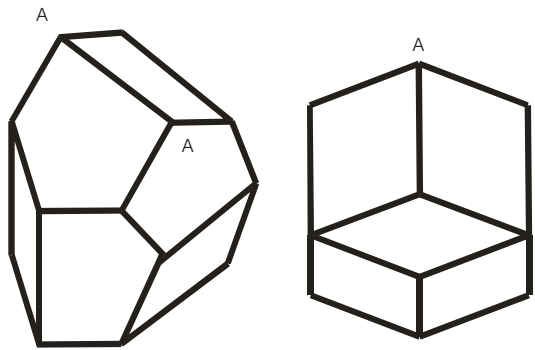

**Figure 4 transformation to the ideal rhombic dodecahedron**

By placing the rhombic dodecahedron in the hyperbolic 3 ball with its six 4-vertices on the sphere at infinity it obtains a hyperbolic metric. Locally, each 3-vertex and edge reverts to having a Euclidean cube geometry so that no geometric singularities occur on identification. Now the rings and thus the singularities of the foliation are pushed out to infinity, and are therefore no longer in the knot compliment. So we have foliated the hyperbolic knot compliment. However we can no longer expect the leaves to be minimal surfaces.

We can consider further links by using more than one cube, or by subdividing the cube into more surfaces. For example if each face is subdivided into four, with a pattern as in figure 5.

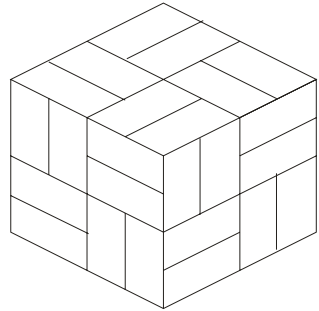

**Figure 5: More elaborate links**

**2.3 Foliating a trivial knot compliment neighborhood**

In Euclidean 3-space, take the four unit cubes with one vertex at the origin and 3 faces in the planes x=0, y=0 and z=0, and lying in the region y≤0. For y>0 the cubes are identified on faces where they touch in space. On the y=0 plane we can place an identification pattern on the cubes. This is shown in figure 6 with the cubes slightly separated.

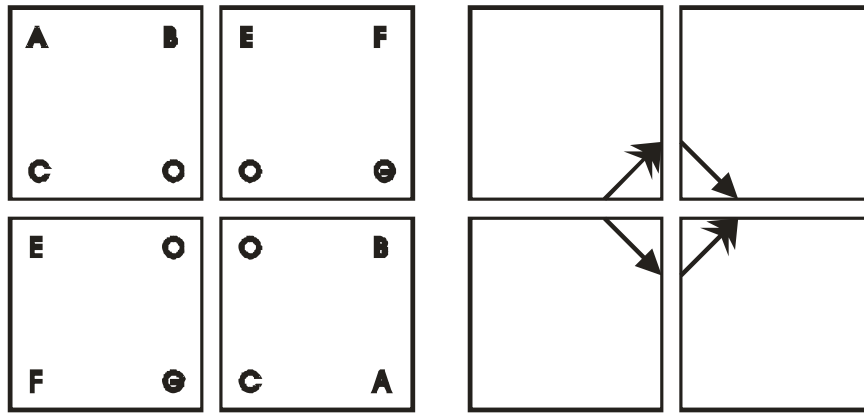

**Figure 6**

The identification pattern given by the letters in the four cubes on the right gives rise to the creation of a non-manifold point at the origin. The link of the origin is a torus not a sphere. Consider removing a pyramid neighborhood of the origin by slightly truncating each cube at the origin. Then the link can be seen as four triangles that all meet at a point on the y axis, behind the origin as we are looking at it in figure 6, and with edges in the plane y=0 having the identification pattern shown. Thus the link of the origin is a torus, not a sphere. If the origin is removed then the four cubes give a neighborhood of the trivial knot in its compliment.

If we ignore the origin which is not a manifold point any longer, this space can be foliated by any direction of foliation (an $\mathbb{RP}^2$ of directions). This is because each isometry in $\mathbb{R}^3$ to make each identification, and thus each holonomy in the identified space, is either the identity for y>0, or a 180° rotation about one of the lines y=0, z=x or y=0, z=-x. The main theorem of this paper concerns constructing cones of general topological surfaces not just tori. First we need to discuss a definition of foliation at the cone points.

## 3   Defining pseudo-foliations and singular pseudo-foliations

We will be concerned with foliation structures whose leaves are area minimizing within local neighborhoods and which are hyperplanes on Euclidean neighborhoods. The spaces we will foliate have geometric singularities in the form of distributional Gaussian curvature and topological singularities in the form of non-manifolds points, such as the cone of $\mathbb{RP}^2$.

We must therefore extend the general definition of foliation to cover interior points in the foliated space which are not manifold points. Before stating the definition we will illustrate how it describes neighborhoods in pseudo-foliations. See figure 7.

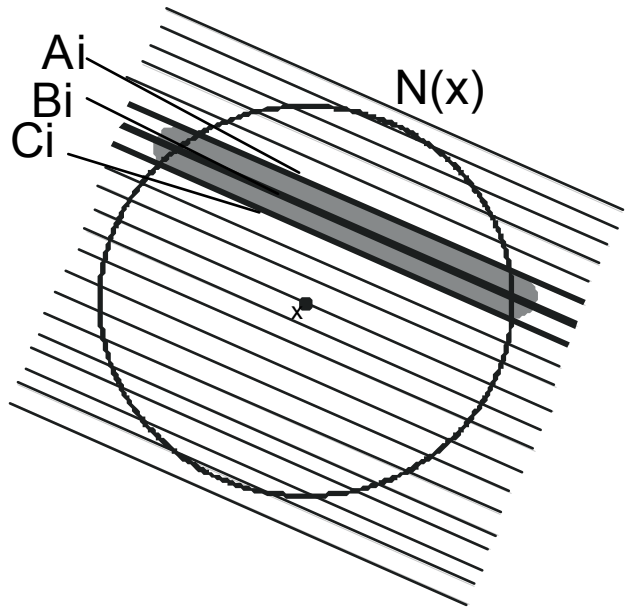

**Figure 7: Definition of pseudo-foliation**

A point x and its neighborhood N(x) are shown together with some of the leaves of the pseudo-foliation and, shaded, a tubular neighborhood $T_t(B_i)$, where B is a leaf and $B_i$ is a component of $(B \cap N(x))$. Note also that for $A_i \in (A \cap N(x))$ and $C_i \in (C \cap N(x))$ $N(x) \cap \partial(T_t(B_i)) = A_i \cup C_i$. We will now state the definition of pseudo-foliation.

**Definition:** *A codimension 1 <u>pseudo-foliation</u> of a space X is a set of submanifolds of X called leaves satisfying the following conditions:*

*(i) Each interior point in X intersects a leaf, and no two leaves intersect.*

*(ii) Each interior point in X has a compact neighborhood N(x) where:*

*a) if $L_i$ is a connected component of $L \cap N(x)$ where L is a leaf of the foliation, then $L_i$ has a one parameter family of tubular neighborhoods $T_t(L_i)$ parameterized by t so that $\bigcup_t T_t(L_i) = L_i$. Also each component of $N(x) \cap \partial(T_t(L_i))$ is a component of $L' \cap N(x)$ where L' is a leaf.*

*b) each component $L_i$ of $L \cap N(x)$ is a component of a $N(x) \cap \partial(T_t(L'_i))$ for some t and some nearby component $L'_i$ of $L' \cap N(x)$ where L' is a leaf.*

*(iii) If $x \in \partial X$ then the pseudo foliation structure at x satisfies the standard definition of a foliation in a neighborhood of a boundary point.*

**Definition:** *A <u>singular pseudo-foliation</u> of a space X is a foliation whereby regular points of a foliation in X have an open neighborhood where the above conditions hold, and on singular sets conditions (i) and (ii) are relaxed to allow leaves to intersect.*

Note that this definition allows for the singularities created when a surface of genus 2 or more is foliated by lines as required by the Poincaré-Hopf index theorem.

## 4      Basic results

### 4.1 Compact surfaces

**Proposition 1:** *It is well known that compact surfaces of non-positive Euler characteristic can be given a piecewise Euclidean metric with one singular point of negative curvature equal to an integer multiple of $2\pi$.*

**Proof:** This can be verified by taking the regular *2n*-gon and identifying its edges to make the desired compact surface with one vertex class. By Gauss-Bonnet, the solid angle will be $2\pi$ plus an angle deficit of $2\pi$ times the Euler characteristic of the surface. Hence and integer multiple of $2\pi$. QED.

We can contrast this special case of negative curvature where the solid angle is an integer multiple of $2\pi$ with the general case. In the case where we have an integer multiple of $2\pi$ we can partition a neighborhood of the vertex into half planes. These can each be foliated by geodesic leaves parallel to their boundaries. This partitioning into half planes also demonstrates by construction that cone points with a piecewise Euclidean neighborhood that can be foliated by geodesics will have a cone angle of an integer multiple of $\pi$.

**Theorem 2:** *Let S be a compact orientable surface. The piecewise Euclidean structure in proposition 1 enables foliations by straight lines on the surface in any direction.*

*There is one singular point of the foliation at the single identified vertex where leaves intersect for surfaces of genus greater than 1. There are no singularities for genus one tori.*

**Proof**: Start with a foliation of a regular 2n-gon in the plane by parallel lines. Pairs of opposite sides of the 2n-gon are identified (by translations). The foliation extends over the edges because parallel edges are identified, and extends around the vertex with a singularity of the required index (1) because the cone angle is an integer multiple of $2\pi$. This can also be verified by examining how foliations of neighborhoods of vertices in the 2n-gon match up as the vertices are identified. The neighborhoods become connected in a cycle and at each transition the foliations match up. QED.

### 4.2 $\mathbb{S}^1$ bundles of surfaces

### 4.2.1 Foliating products of $\mathbb{S}^1$ with compact surfaces

**Theorem 3:** *We can foliate products of $\mathbb{S}^1$ with compact surfaces.*

**Proof:** Consider the product of $\mathbb{S}^1$ with a topological surface. One foliation is just to take each copy of the surface as a leaf. Taking the metric to be the product metric this foliation will be minimal with no singular points. We also know, above, that we can foliate an orientable surface with straight lines and in the product this will give leaves that are products of $\mathbb{S}^1$ and the straight line leaves. These foliations will, in general, be singular on a set $\mathbb{S}^1 \times \{p\}$, where p is the cone point. Leaves of the foliation will intersect on this set. QED.

**Theorem 4:** *When the cone angle in a cone is an integer multiple of $2\pi$ the product of the cone with an interval will admit a foliation by hyperplanes in any direction.*

**Proof:** This can be seen by taking n copies of a unit cube in $\mathbb{R}^3$, where n is the integer multiple of $2\pi$ giving the curvature. Put the same foliation on each copy. Then cut each cube in the same half plane passing through the three midpoints of the top face the bottom face and the same side face. Now identify the cubes in a cycle creating the desired product structure. To do this for the example of three cubes, see figure 8 (see also section 1.3).

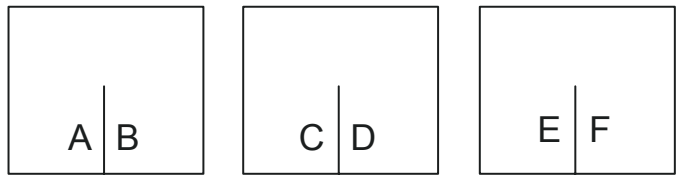

**Figure 8: End view of three cubes to be identified**

Three cubes are shown end on with the line to the center indicating the removed rectangle intersecting the square face shown. So A, B, C, D, E and F each represent a side of a removed rectangle. The identity map would identify A with B, C with D and E with F. It would also, in the process, replace the missing rectangles. To connect the cubes we will add three rectangles, but this time with a different identification pattern. A will be identified with D, as one rectangle is added, C with F as an other rectangle is added and E with B for the final rectangle. This process can be extended to n cubes. QED.

**Theorem 5:** *Each leaf of the foliation in theorem 4 is an area minimizer.*

**Proof:** A calibration form orthogonal to the leaves extends across the new product structure. The singular set has area measure zero, so can be ignored. The calibration form proves that all leaves are minimizers. QED.

**Theorem 6:** *The space of piecewise planar foliations of the above products of surfaces with $\mathbb{S}^1$ is $\mathbb{RP}^2$, a full Grassmannian worth.*

**Proof:** Take the polygon cross the unit interval as a subset of $\mathbb{R}^3$. Choose any unit vector in $\mathbb{R}^3$ and use its normal planes to foliate the polygon cross the interval. The space of such piecewise planar foliations is $\mathbb{RP}^2$. The identifications on faces and the edges match up as in theorem 4. QED.

### 4.2.2 Other $\mathbb{S}^1$ bundles.

If we take a regular 2n-gon, then we can apply a 180° rotation about the center. This will preserve all foliations when it becomes the holonomy map in the $\mathbb{S}^1$ bundle. For a circle bundle over a torus, the monodromy map can be a translation. The space of piecewise planar foliations of all these spaces and their connect sums, see below, is also $\mathbb{RP}^2$, a full Grassmannian worth.

We can generalize the comparison made between the annulus and moebuis band under the locally piecewise flat condition (section 2.2). If the bundle is twisted, then there will be a foliation because the holonomy around $\mathbb{S}^1$ will be a reflection. If however the holonomy is a rotation, then a foliation can only exist if the rotation angle is zero or π.

### 4.3 Extending foliations across connect sum operations.

**Theorem 7:** *Suppose 3-manifolds M and N each have a piecewise Euclidean structure. Also assume that M has a foliation by area minimizing surfaces which on Euclidean neighborhoods, as subsets of $\mathbb{R}^3$, have leaves which are parallel hyperplanes. Then a connect sum can be constructed between M and N so that the foliation M extends uniquely to M#N. This foliation on M#N will also have the property of having leaves which are area minimizing surfaces which on Euclidean* *neighborhoods, as subsets of $\mathbb{R}^3$, are parallel hyperplanes. This process also works for other dimensions.*

**Proof:** We adapt the construction of the connect sum, from the usual procedure of removing a 3- ball from each and identifying the 2-sphere boundaries, by using the process in 1.3. We remove two planar discs from each 3-manifold and compactify as in 1.3 to create a topological sphere boundary in each manifold. Then proceed as usual with the identifications on those spheres. Geometrically these spheres are immersed in planar discs.

A minimal foliation and its calibration form can be extended across the identification by finding two corresponding directions on each side of the connect sum. The correspondence of directions ensures that if the foliation is deformed near a disc, any reduction of area on one of the components will be compensated for by an equivalent increase on the other. See figure 9. The diagonal leaf switches from the left component above the thick line to the right component below the thick line (representing the discs) with the identifications as shown by the arrows. Moving the contact point with the thick line from above to the left will shorten the leaf above the line but will move the leaf to the right and lengthen it below the line.

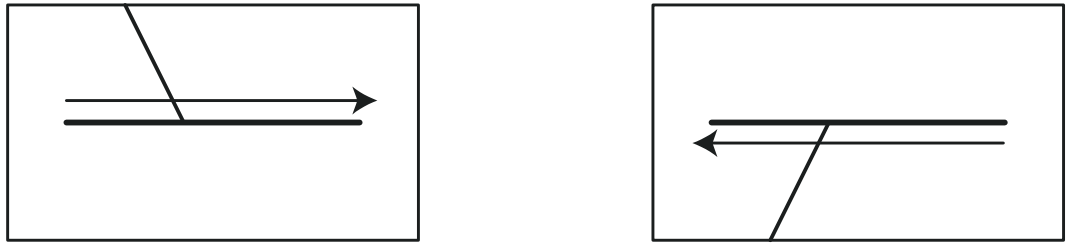

**Figure 9: Extending a foliation across the connect sum operation**

This same process can be carried out for *n*-manifolds where *n* is greater than 1. The identification will be carried out on n-spheres immersed in embedded hyperplane *(n-1)*-balls. QED.

## 5 Construction of piecewise Euclidean structures for cones of compact connected surfaces that admit pseudo-foliations.

### 5.1 No trivial solutions

If we take a compact surface S embedded in $\mathbb{R}^n$, there is a trivial cone structure in $\mathbb{R}^{n+1}$ which is given by the set $K_S=\{(ts,t) \subset \mathbb{R}^n \times \mathbb{R}: s \in S \subset \mathbb{R}^n, t \in [0,\infty)\}$. $K_S$ can inherit the Euclidean metric from $\mathbb{R}^{n+1}$. There are two types of foliation we can consider. One is to take level sets of t, where each leaf an image of S under some homothety. This gives a foliation away from t=0, but at t=0, there is an isolated point not part of a leaf. Also this foliation is not minimal everywhere for any value of t, as the vector field –x, will reduce area on every leaf linearly with time. So the leaves are not area stationary.

The other trivial type of foliation of the cone would be to take a foliation on S and simply take the cone of each leaf to give the leaves of the foliation of $K_S$. This is a foliation away from t=0, and the leaves may even be minimal surfaces, however at t=0 all the leaves intersect.

**5.2 The construction to give foliations by minimal surfaces**

**Theorem 8:** *A cone of any compact connected boundary-less surface can be given a piecewise Euclidean structure that admits a pseudo-foliation by piecewise planar area minimizing surfaces.*

**Proof:** We shall now give a general construction for the cone of any compact surface with one connected component using the method of section 1.3. To start with a simple example, take a cube as in the upper left of figure 10 equipped with a Euclidean metric structure. Remove triangle LMN, where L is in the interior and M and N are on the bottom boundary of the cube.

The next stage is to place identifications on the two copies of the triangle so that the new boundary is removed. The identifications are made according to the arrows labeled 'a' and 'b' at each horizontal level, but not between horizontal levels. Note that as the triangle gets thinner nearer L the lengths of 'a' and 'b' get shorter, but the identification pattern remains the same. This creates the cone of a surface with L being the cone point.

Figure 10 shows the cone of a sphere (bigon with boundary identification $aa^{-1}$) on the left and the cone of a projective plane (bigon with boundary identification aa) on the right. To see the projective plane look at the entire boundary of the cube. Also below in figure 10 is the cone of a torus (rectangle with boundary identification $aba^{-1}b^{-1}$).

One pseudo-foliation is simply where all horizontal planes in the cube are leaves. The proof that this is minimizing involves a calibration argument using the form *dz*, where *x*, *y* and *z* are orthogonal directions with *x* and *y* horizontal and *z* vertical. See introduction for a full explanation of calibration and free boundary. Locally at every point away from triangle LMN the calibration form *dz* proves that the horizontal surfaces are minimal. If we treat the triangle LMN as free boundary, see figures 10 and 11, then we can use the calibration form *dz* to prove minimality of the leaves in that free boundary case. The horizontal leaf shown in figure 11 has boundary shown in thick lines. The outside is fixed boundary and the thick line inside is where the leaf touches the free boundary formed by triangle LMN. When we say the leaf is a minimizer by calibration we mean that any surface with the same fixed boundary and free boundary will have areas greater than or equal to our leaf.

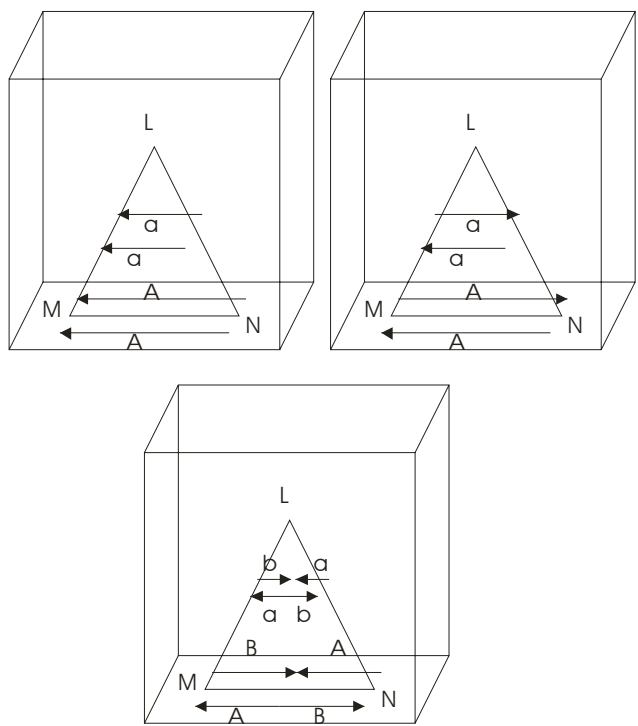

**Figure 10: Cones of the sphere, $\mathbb{RP}^2$ and a torus**

Remark: There is a change in topology of the leaves of the above pseudo-foliations as the horizontal leaves travel down over the non-manifold point from above to below. A handle or cross cap is added to the surface.

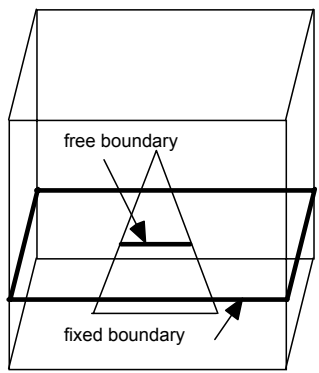

**Figure 11, A leaf with fixed and free boundary.**

We now apply the short cut lemma (lemma 9) below to prove theorem 8 by saying that as our leaves are minimizers after free boundary has been added, then they are minimizers in the original cone space.

**5.3 The short cut lemma**

We now need a definition of adding free boundary so we may introduce the short cut lemma. Say we have a length, area or volume minimizing problem for a submanifold with fixed boundary in an ambient manifold. We can introduce a set into the ambient space on which our minimizer can also have boundary. We shall call any boundary components of this set free boundary components. For example see figure 12. Two pictures are shown of two boats, A and B. A swimmer has to get from one boat to the other and there is no current in the water. On the left there is only clear water to swim through in a straight line. Thus A and B are the fixed boundaries of the path. Once the island is introduced, the swimmer can swim to one side of the island, walk to the other side and swim the rest of the way. The island perimeter is a free boundary for the swimmer as the swimmer is free to choose where to land and where leave the island from. Altogether, in this case, the swimming sections have four boundary points, A, B, and the two points on the island.

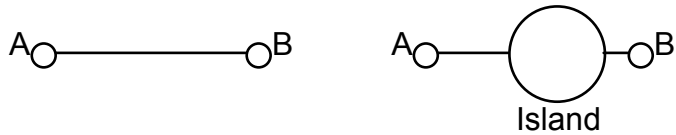

**Figure 12**

We shall consider a second example which we can regard as a one dimensional prototype for the subsequent lemma.

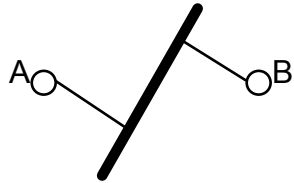

**Figure 13**

Say instead of having a circular island, the swimmer can land and walk along a pier between the boats shown in bold (Figure 13). This pier represented by the thick line is added free boundary, although it has no interior points as we are regarding it as a line. Now we can state and prove the lemma.

**Lemma 9:** *If a surface S with fixed boundary in a space X becomes a minimizer in a space Y obtained by defining a subset of X to become free boundary, then it was already a minimizer in the original space X.*

For example see figure 11. The space X was the space before the triangle was determined to be free boundary. The space Y is the space as shown in figure 11 with the triangle defined as free boundary. The horizontal leaf is a minimizer in the space Y by calibration, and hence also in X by the short cut lemma.

**Proof:** If a surface S is an area minimizer among candidates for a space Y, and all candidate surfaces for space X are candidates for space Y, and the area forms in Y and X are equal up to a set of measure zero, then S is a minimizer in space X. QED

By proving lemma 9, we have completed the proof of theorem 8.

## 5.4 Searching for the space of our pseudo-foliations by piecewise planar area minimizing surfaces

Using the short cut lemma we can also prove that two other piecewise planar pseudo-foliations, possibly singular pseudo-foliations, are area minimizing.

**Theorem 10:** *There are three piecewise planar pseudo-foliations, singular pseudo-foliations, of the constructed space in theorem 8, for the cone of every connected compact surface, whose leaves are mutually orthogonal on every Euclidean neighborhood.*

**Proof:** We will use the coordinate axes in figure 14 to identify our three piecewise planar, possibly singular, pseudo-foliations. This gives us, possibly singular, pseudo-foliations by surfaces with equations of the form $x=k$, by surfaces with equations of the form $y=k$, and by surfaces with equations of the form $z=k$. For the surfaces of the form $x=k$ and $z=k$, the triangle LMN is made a free boundary, and then the short cut lemma is used. For the $y=k$ surfaces an orthogonal set of free boundaries is used. All these surfaces contain lines parallel to the y axis. These intersect LMN at points where arrows such as 'a' and 'b' (see figure 10) begin and end. This effectively cuts the space into purely Euclidean piecewise flat sections in which the short cut lemma is applied. QED.

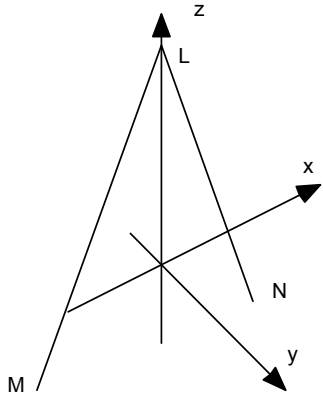

**Figure 14: The positions of the coordinate axes with respect to triangle LMN.**

Note that for equations $x=k$ and $y=k$ the foliation structure is a singular pseudo-foliation on the cone of the heads and tails of identification arrows for cones of surfaces of genus 0 or higher. This corresponds to the way in which vertices are identified in polygons to make surfaces.

**Theorem 11:** *The space of area minimizing piecewise planar foliations, possibly singular, near a cone of a compact connected surface, constructed as in theorem 8 contains a circle union a disjoint point.*

**Proof:** We can also show with a different use of calibration [Morgan], [Harvey and Lawson] that equations of the form $ay + bz = k$, $a\neq 0 \neq b$, describe minimizing surfaces for the constructed cones. For this we need to show that deformations, S, of the surface T as shown in figure 15 do not decrease area. The deformation shown pushes the surface up on the left in front of the triangle LMN and correspondingly pushes it up on the right at the rear. This surface does not have boundary on the triangle LMN, because of the identification pattern.

To use the short cut lemma, we convert triangle LMN to free boundary and remove the identification pattern. This creates new boundary in the form of two semi-circles shown in figure 15.

A calibration form to prove the plane is a minimizer can be represented as -*adz* + *bdy*. To use calibration we must remove the new boundary by adding surface $d_1 \cup d_2$, the black semi-circles shown. Fortunately when we integrate (in (*) below) using the calibration form, the contribution of the black semi-circles cancel as they have opposite orientation (note that these black semi circles will not always be perpendicular to the surface).

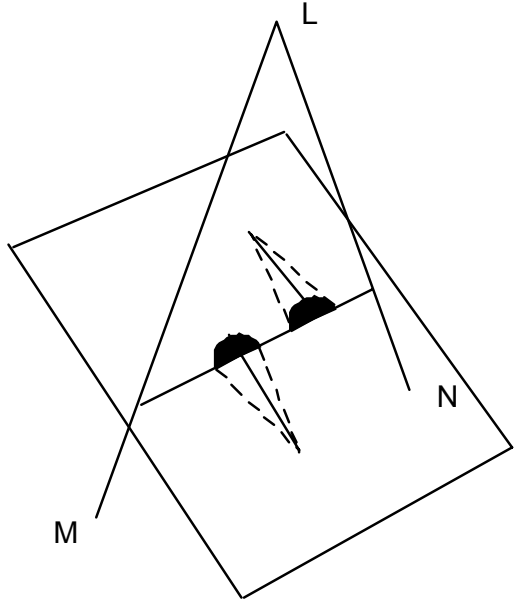

**Figure 15: Deformation of plane ay + bz = k**

If T is the rectangular portion of the calibrated plane *ay* + *bz* = *k*, as shown, and S is the deformation of T union the two semicircles $d_1$ and $d_2$, we can write:

$$Area(S) \underset{\substack{Form \\ doesn't \\ generally \\ calibrate \\ S}}{\geq} \int_S -adz + bdy \underset{\substack{Semicircles' \\ orientations \\ opposite}}{=} \int_{S \cup d_1 \cup d_2} -adz + bdy \underset{\substack{Stoke's \\ Theorem}}{=} \int_T -adz + bdy \underset{\substack{Form \\ calibrates \\ T}}{=} Area(T)$$

(*)

We have shown that the deformed surface S will not have less area than the original T. Hence the original surface T is a minimizer with the free boundary. Finally the short cut lemma implies that T is also a minimizer in the cone of a surface.

In the Grassmann bundle we now have a circle of values for the directions of minimizing foliation planes with equations *ay* + *bz* = *k*, union a point for the planes *x* = *k*. This proves that the space of area minimizing piecewise planar pseudo-foliations near a cone of a compact connected surface contains a circle union a point. QED.

## Acknowledgements

Thanks to Robert Hardt, John Hempel, Thierry De Pauw, Robin Forman, Michael Wolf, Al Marden, Frank Morgan, Richard Evans and Shelley Harvey. Particular thanks to Penny Smith for suggesting the study of hyperbolic 3 manifolds.

## References

Banchoff, T. Foliations of knot compliments in the bicylinder boundary, *Bol. Soc. Brazil. Mat*. 5 (1974), 31-43

Dippolito, P. Codimension one foliations of closed manifolds. *Ann. Of Math*. 107,(1978), 403-453.

Epstein, David and Gunn, Charlie. *Not knot*. Video and Booklet. The Geometry Center University of Minnesota. Distributed by A.K. Peters, Wellesley Massachusetts.

Gabai, David 3 lectures on foliations and laminations on 3-manifolds. Laminations and foliations in dynamics, geometry and topology (Stony Brook,NY, 1998), 87--109, *Contemp. Math.,* 269, Amer. Math. Soc., Providence, RI, 2001.

Gabai, David, William H. Kazez Group negative curvature for 3-manifolds with genuine laminations. *Geom. Topol.* 2 (1998), 65-77.

Harvey, Reese, and Lawson, H. Blaine, Jr. Calibrated geometries. *Acta Math.* 148 (1982), 47-157

Jones, Kerry N. Geometric structures on branched covers over universal links. Geometric topology (Haifa,1992), 47-58, *Contemp. Math.,* 164, Amer. Math. Soc., Providence, RI, 1994.

Meeks, William, III; Simon, Leon; Yau, Shing Tung. *Embedded minimal surfaces, exotic spheres, and manifolds with positive Ricci curvature.* Ann. of Math. (2) 116 (1982), no. 3, 621-659.

Morgan, Frank. *Geometric measure theory: a beginner's guide*. Academic Press, 2000.

Munkres, James R. *Elements of Algebraic Topology*. Addison Wesley, 1984.

Jaco, William(1-OKS); Rubinstein, J. Hyam. *0-efficient triangulations of 3-manifolds*. J. Differential Geom. 65 (2003), no. 1, 61--168.

Stocking, Michelle. *Almost normal surfaces in 3-manifolds*. Trans. Amer. Math. Soc. 352 (2000), no. 1, 171-207.

Sullivan, Dennis. A homological characterization of foliations consisting of minimal surfaces. Insititute des Hautes Études Scientifiques. *IHES/M/78/202 Comment. Math. Helv*. 54 (1979), 218-223.

Thurston, William, P. *The geometry and topology of 3-manifolds*. Unpublished notes.